\documentclass{article}

\usepackage{arxiv}

\usepackage[utf8]{inputenc} 
\usepackage[T1]{fontenc}    
\usepackage{hyperref}       
\usepackage{url}            
\usepackage{booktabs}       
\usepackage{amsfonts}       
\usepackage{nicefrac}       
\usepackage{microtype}      
\usepackage{lipsum}		
\usepackage{graphicx}
\usepackage{amsmath}
\usepackage{subfigure}

\begin{document}

\title{The fundamental solutions of the curve shortening problem via the Schwarz function}


\author{ N.R. McDonald \\
	Department of Mathematics\\
	University College London\\
	London, WC1E 6BT \\
	\texttt{n.r.mcdonald@ucl.ac.uk} \\
}


\renewcommand{\headeright}{}
\renewcommand{\undertitle}{}

\hypersetup{
pdftitle={Curve shortening solutions and the Schwarz function},
pdfauthor={N.R. McDonald},
pdfkeywords={curve shortening, Schwarz function},
}

\title{The fundamental solutions of the curve shortening problem via the Schwarz function
}


\author{ N.R. McDonald \\
	Department of Mathematics\\
	University College London\\
	London, WC1E 6BT \\
	\texttt{n.r.mcdonald@ucl.ac.uk} \\
}




\maketitle

\begin{abstract}
Curve shortening in the $z$-plane in which, at a given point on the curve, the normal velocity of the curve is equal to the curvature, is shown to satisfy $S_tS_z=S_{zz}$, where $S(z,t)$ is the Schwarz function of the curve. This equation is shown to have a parametric solution from which the known explicit solutions for curve shortening flow; the circle, grim reaper, paperclip and hairclip, can be recovered.

\end{abstract}

\section{Introduction}
\label{intro}
The evolution of a smooth curve ${\pmb\gamma}(\textbf{x},t)$ in the plane with given initial shape ${\pmb\gamma}_0$  in which the normal velocity at a given point on the curve is proportional to the curvature $\kappa$ at that point, is known as the curve shortening problem:
\begin{align}
\frac{\partial{\pmb\gamma}}{\partial t}&=\kappa{\textbf{n}(\textbf{x},t)},\nonumber \\
{\pmb\gamma}({\textbf{x}},0)&={\pmb\gamma}_0,
\label{csivp}
\end{align}
where $\textbf{n}$ is the normal.
The problem, and its higher dimensional generalisation, has attracted much attention since the 1980s. Important features of the curve evolution, such as the shrinking of any closed, embedded curve first to a convex curve and then to a round point in finite time, have been established e.g. \cite{gage,grayson}. Curve shortening (\ref{csivp}), and its variants, also has practical application; for example in the late-time evolution of Hele-Shaw free boundary flow in the presence of surface tension \cite{dall}. In the time reversed sense when the curve lengthens, the solutions have relevance to viscous fingering and crystal growth e.g. \cite{nak} where the connection to Saffman-Taylor fingering is made.

There are only four known, explicit, time-dependent solutions describing curve shortening: the bounded, closed curve solutions of the (i) circle and (ii) paperclip; and the unbounded, (iii) steadily translating grim reaper, and the (iv) hairclip solution--see e.g. \cite{nak,broad,tsai}. A straight line is a trivial solution and remains stationary. More exotic, but not explicit, solutions taking the form of spirals which rotate and/or expand are described in \cite{hall}.

\begin{figure}
\centering
\mbox{\subfigure{\includegraphics[width=2.8in]{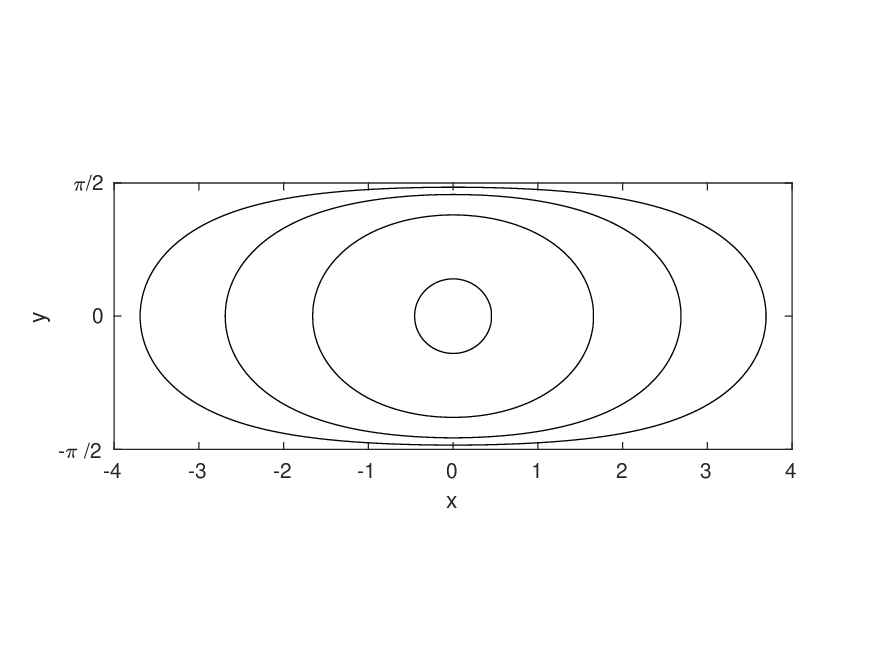}}\quad
\subfigure{\includegraphics[width=2.8in]{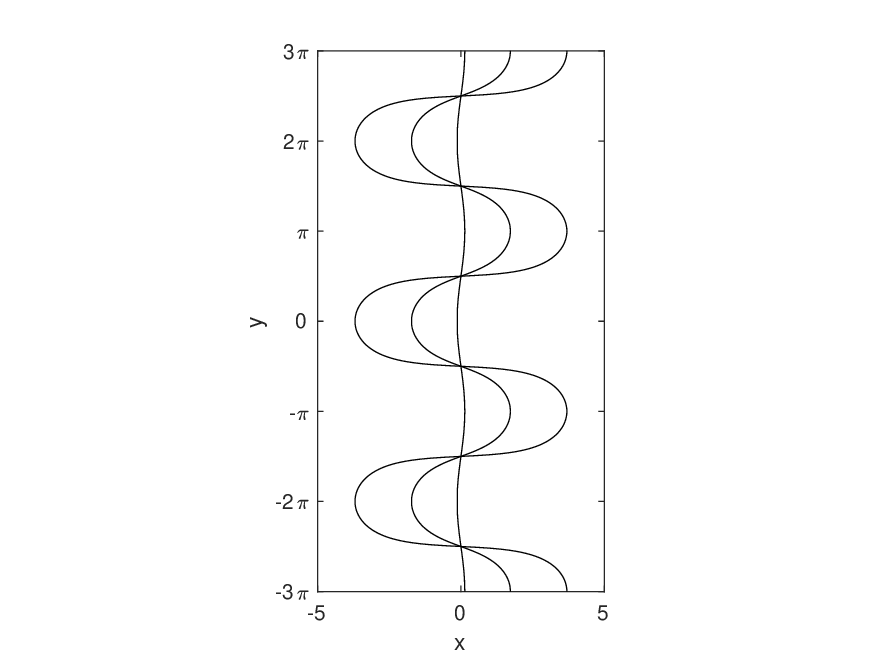} }}
\caption{The evolution of the paperclip solution (left) collapsing toward a circle shown at times $t=-3,-2,-1,-0.1$. The evolution of the hairclip solution (right) collapsing toward the vertical line $x=0$ shown at times $t=-3,-1,2$.} \label{fig1}
\end{figure}

Solutions (ii)-(iv) can be written (see e.g. \cite{tsai} after a suitable rotation of coordinates and choice of timescale)
\begin{itemize}
\item[] Grim reaper: $x=t-\log(\cos y)$;
\item[] Paperclip: $x= \cosh^{-1}(\exp(-t)\cos y)$; 
\item[] Hairclip: $x=- \sinh^{-1}(\exp(-t)\cos y)$.
\end{itemize}
Sketches of evolving paperclip and hairclip solutions are shown in Fig. \ref{fig1}.

There is interest in deriving these exact solutions systematically using simple methods. One approach has been to show the fundamental solutions can be obtained from a nonlinear diffusion equation \cite{nak,broad} which can be solved by functional separation of variables in conjunction with symmetry assumptions. More recently, \cite{tsai} presents two methods for deriving the four fundamental solutions. The first is based on seeking explicit solutions to (\ref{csivp}) in which the function being sought is itself a function satisfying the classic  heat equation $h_t=h_{xx}$. Assuming $h_t$ is either constant or linear in $h$ gives each of the four fundamental solutions. The other approach detailed in \cite{tsai} is an implicit method in which the solution is sought in the form of separable functions of the three variables $x,y,t$.

An alternative and systematic procedure for obtaining the exact solutions for time-dependent curve shortening is presented here. The key idea is to represent the curve in the complex plane using the Schwarz function \cite{davis}. Then, in Sec. \ref{sec:2}, well-established results for representing the geometric properties of the curve in terms of the Schwarz function are used to derive a PDE governing the evolution of the curve.  Sec. \ref{sec:3a} shows that the self-similar, collapsing circle is an exact solution of the PDE, and Sec. \ref{sec:3b} uses elementary methods to derive a steadily translating parametric solution of the PDE which is equivalent to the grim reaper solution. In Sec. \ref{sec:3c} and \ref{sec:3d}, the parametric solution is generalised and shown to yield the paperclip and hairclip solutions. Sec. \ref{sec:4} reconsiders the steadily translating problem, and derives directly the grim reaper solution from the Schwarz function equation cast in a  complex plane moving with the solution curve. The same approach is used to find a new exact solution taking the form of a steadily translating parabola to the problem in which the normal velocity of the curve is given by the cube root of the curvature.

\section{The governing equation in terms of the Schwarz function}
\label{sec:2}
The Schwarz function $S(z)$ of an analytic curve $\pmb\gamma$ in the $z$-plane is the unique function which is analytic in the neighbourhood of $\pmb\gamma$ and such that $S(z)={\bar z}$ on $\pmb\gamma$ \cite{davis}. For example, a circle of radius $a$, centred at the origin, has $S(z)=a^2/z$. Representing curves using $S(z)$ has proved remarkably successful in obtaining exact solutions in fluid mechanical applications where the evolving curve may represent, for example, the interface between fluids of different viscosities such as occurs in a Hele-Shaw cell (e.g. \cite{gust,howi2,mineev}), or, in 2D vortex dynamics, where the interface separates flows with different constant vorticities (e.g.\cite{saff,crowd,ricc}). The common theme in these applications is that the dynamics of a planar curve can be cast in terms of an equation satisfied by the Schwarz function which, when solved, determines the shape of the curve. Often the success of this approach is owing to the fact that $S(z)$ is an analytic function, at least in the neighbourhood of $\pmb\gamma$, enabling methods of complex analysis, such as conformal mapping, to be used. 

In terms of $S(z)$, the curvature $\kappa$ of $\pmb\gamma$ is \cite{davis}
\begin{equation}
\kappa=\frac{i}{2}\frac{S_{zz}}{(S_z)^{3/2}},
\label{curv}
\end{equation}
where subscripts denote partial derivatives.
For an evolving curve $S=S(z,t)$, the normal velocity $v_n$ at point $z$ on $\pmb\gamma$ is (e.g. \cite{mineev})
\begin{equation}
v_n=\frac{i}{2}\frac{S_{t}}{(S_z)^{1/2}}.
\label{norm}
\end{equation}
Curve shortening equates (\ref{curv}) and (\ref{norm}), giving a nonlinear PDE for $S(z,t)$:
\begin{equation}
S_tS_z=S_{zz}.
\label{schwarz}
\end{equation}
The curve shortening problem involves solving (\ref{schwarz}) subject to the initial shape of the curve $S(z,0)$.

An immediate feature of (\ref{schwarz}) is, given ${\bar z}=S(z)$, that it is invariant to the scaling $z\to\lambda z$ ($\lambda\in \Re$), which implies $\bar{z}=S\to\lambda S$,  and $t\to\lambda^2 t$ -- see also \cite{tsai}. Moreover, as expected, the governing equation (\ref{schwarz}) is invariant to rotations through angle $\sigma$, since 
$z\to e^{i\sigma}z$ implies $S(z)\to e^{-i\sigma}S$, rendering (\ref{schwarz}) invariant.

Letting $s$ be the arclength parameter of $\pmb\gamma$, it can be shown (e.g. \cite{davis}) $dz/ds=1/\sqrt{S_z}$, and upon conjugating $dS/ds=\sqrt{S_z}$. These relations can be used in (\ref{schwarz}) to establish the known connection between the curve shortening flow and the heat equation (e.g. \cite{tsai}), which in term of the Schwarz function is
\begin{equation}
S_t=2S_{ss}.
\label{heat}
\end{equation}

\section{Explicit solutions of the curve shortening flow}
\label{sec:3}

\subsection{The collapsing circle}{\label{sec:3a}}
By symmetry, an initially circular curve remains circular under the action of curve shortening: this is a known exact (self-similar) solution. It is simple to demonstrate using (\ref{schwarz}): a circle of radius $a(t)$ with (without loss of generality) centre at the origin, has Schwarz function $S(z,t)=a(t)^2/z$ which satisfies (\ref{schwarz}) exactly provided $a{\dot a}=-1$, where $\dot a$ is the time derivative of $a(t)$. This is equivalent to ${\dot A}=-2\pi$ where $A=\pi a^2$ is the area enclosed by the circle, and consistent with the well-known result that the area enclosed by an arbitrary closed
plane curve
decreases at the rate of $2\pi$ per unit time e.g. \cite{grayson}.

\subsection{The grim reaper}{\label{sec:3b}}
Let a curve translate steadily, with unit velocity, in the positive real direction, and write it as
\begin{equation}
z=t+f(\zeta),
\label{translate}
\end{equation}
where $\zeta=\exp{(i\theta)}$, $-\pi\le\theta\le \pi$, and $\bar{f}=f$. The latter condition implies the curve is symmetric about the real axis. Now 
\begin{equation}
S(z,t)=\bar{z}=t+\overline{f({\zeta})}=t+f(\zeta^{-1}).
\label{grs}
\end{equation}
It follows from (\ref{translate}) and (\ref{grs}) that
\begin{equation}
S_t=1+\frac{f'(\zeta^{-1})}{\zeta^2f'(\zeta)},\quad \textrm{and}\quad S_z=-\frac{f'(\zeta^{-1})}{\zeta^2f'(\zeta)}.
\label{firstderivs}
\end{equation}
Also,
\begin{align}
S_{zz}&=\frac{1}{f'(\zeta)}\partial_\zeta S_z,\nonumber\\
&=\frac{1}{f'(\zeta)}\left ( \frac{f''(\zeta^{-1})}{\zeta^4f'(\zeta)}+\frac{2f'(\zeta^{-1})}{\zeta^3f'(\zeta)}+\frac{f'(\zeta^{-1})f''(\zeta)}{\zeta^2f'(\zeta)^2}\right ).
\label{secderiv}
\end{align}

Substituting (\ref{firstderivs}) and (\ref{secderiv}) into (\ref{schwarz}) and simplifying gives
\begin{equation}
\frac{d}{d\zeta}\log\frac{f'(\zeta^{-1})}{f'(\zeta)}=\frac{2}{\zeta}+f'(\zeta)+\frac{f'(\zeta^{-1})}{\zeta^2}.
\label{feq2}
\end{equation}
Let $h(\theta)=e^{i\theta}f'(e^{i\theta})$ and using $\partial_\zeta=i e^{i\theta}\partial_\theta$, (\ref{feq2}) gives a functional-differential equation for the complex-valued $h(\theta)$:
 \begin{equation}
\frac{d}{d\theta}\log\frac{h(-\theta)}{h(\theta)}=i\left [ h(\theta)+h(-\theta)\right ].
\label{feqh}
\end{equation}
Since ${\bar f}=f$, it follows that $h(-\theta)={\overline{h(\theta)}}$ and that the real part of (\ref{feqh}) is satisfied by inspection with both sides of the equation vanishing.

 Letting $h(\theta)=r(\theta)\exp(i\alpha(\theta))$ and substituting into (\ref{feqh}) gives
\begin{equation}
\frac{d\alpha}{d\theta}=-r(\theta)\cos\alpha(\theta).
\label{alphaeq}
\end{equation}
But $h(\theta)=r(\theta)\exp(i\alpha(\theta))=-\alpha'(\theta)[1+i\tan(\alpha(\theta))]$ and is the general solution to (\ref{feqh}). Now $h(\theta)=e^{i\theta}f'(e^{i\theta})$ implies $df/d\theta=ih$ and integrating gives
\begin{equation}
f(\theta)=-\log\left ( 1+e^{2i\alpha(\theta)}\right ),
\label{general}
\end{equation}
where a constant of integration (which is real since ${\bar f}=f$) has been set to zero since it only serves to shift the curve along the real axis.

Equations (\ref{translate}) and (\ref{general}) gives
\begin{equation}
z=t-\log\left ( 1+e^{2i\alpha(\theta)}\right ),
\label{general2}
\end{equation}
and taking real and imaginary parts of (\ref{general2}) gives $x=t-\log(\cos y)$ which is the grim reaper solution. Note that (\ref{general2}) does not admit embedded curve solutions; the curve must close at infinity. That is, $\alpha(\pm\pi)=\pm \pi/2$. The choice of function $\alpha(\theta)$ is immaterial so long as it is a monotonic function such that $-\pi/2\le\alpha(\theta)\le\pi/2$. Choosing the simple form $\alpha(\theta)=\theta/2$ gives the final parametric form of the grim reaper solution 
\begin{equation}
z=t-\log\left ( 1+e^{i\theta}\right ),\quad -\pi\le\theta\le\pi.
\label{grparam}
\end{equation}

\subsection{The paperclip}{\label{sec:3c}}

The grim reaper solution in parametric form (\ref{grparam}) suggests a general form of solution to (\ref{schwarz}):
\begin{equation}
z=-\log\left ( a(t)+e^{i\theta}\right )+g(t),
\label{pc1}
\end{equation}
where $a(t)>0$ (without loss of generality) and $g(t)$ are real-valued functions to be determined. First, the choice $a(t)>1$ is made and it is shown this corresponds to the paperclip solution. In \ref{sec:3d}, the choice $0< a(t)<1$ is shown to correspond to the hairclip solution.

Seeking a closed curve solution which is symmetric about the imaginary axis implies $z|_{\theta=0}=-z|_{\theta=\pi}$. From (\ref{pc1}), $g(t)=(1/2)\log(a^2-1)$, and the corresponding Schwarz function is
\begin{equation}
S(z,t)=-\log\left ( a(t)\zeta+1\right )+\log\zeta+\frac{1}{2}\log(a^2-1),
\label{pcschwarz}
\end{equation}
where $\zeta=e^{i\theta}$.

The first step is to show from (\ref{schwarz}) that $a^2=1/(1-\exp(2t))$ and so $-\infty <t<0$, corresponding to what is known as an {\it ancient} solution. From (\ref{pcschwarz})
\begin{equation}
S_\zeta=\frac{\zeta_z}{\zeta(a\zeta+1)}=-\frac{\zeta+a}{\zeta(a\zeta+1)},
\label{sz}
\end{equation}
and from (\ref{pc1}) and (\ref{pcschwarz})
\begin{equation}
S_t=\frac{-({\dot a}\zeta^2+{\dot a})(a^2-1)+(a\zeta^2+2\zeta+a)a{\dot a}}{\zeta(a^2-1)(a\zeta-1)}.
\label{dsdt}
\end{equation}
Also,
\begin{align}
S_{zz}&=\zeta_z\frac{d}{d\zeta}\left ( \zeta_z\frac{1}{\zeta(a\zeta+1)}\right),\nonumber\\
&=-\frac{(\zeta+a)(a\zeta^2+2a^2\zeta+a)}{\zeta^2(a\zeta+1)^2}.
\label{seconderivpc}
\end{align}
Substituting (\ref{sz}), (\ref{dsdt}) and (\ref{seconderivpc}) into (\ref{schwarz}) and simplifying gives an ordinary differential equation for $a(t)$:
\begin{equation}
{\dot a}=a(a^2-1),
\label{dadt}
\end{equation}
with solution
\begin{equation}
a(t)^2=\frac{1}{1-k^2e^{2t}},
\label{asol}
\end{equation}
where $k$ is a real constant. The choice of $k$ simply determines the termination time 
of the solution, when $a(t)$ becomes singular. Without loss of generality this is chosen to be $t=0$ and therefore $k=1$. Thus
\begin{equation}
z=-\log\left ( a(t)+e^{i\theta}\right )+\frac{1}{2}\log(a^2-1),
\label{pc2}
\end{equation}
where $a^2=1/(1-\exp(2t))$ and $-\pi\le\theta\le\pi$, is a time-dependent, closed-curve, parametric solution of the curve shortening problem.

The next step is to demonstrate the equivalence of (\ref{pc2}) to  the paperclip solution $\exp(-t)\cos y=\cosh x$. 
The real and imaginary parts of (\ref{pc2}) give
\begin{align}
x&=\log\sqrt{a^2-1}-\log\sqrt{a^2+2a\cos\theta+1},\nonumber\\
y&=-\tan^{-1}\left( \frac{\sin\theta}{a+\cos\theta}\right ).
\label{xandy}
\end{align}
Hence
\begin{align}
e^{-t}\cos y&=e^{-t}\frac{a+\cos\theta}{a^2+2a\cos\theta+1}\nonumber\\
&=\frac{a(a+\cos\theta)}{\sqrt{a^2-1}\sqrt{a^2+2a\cos\theta+1}},
\label{lhs}
\end{align}
and
\begin{align}
\cosh x&=\frac{1}{2}\left ( \frac{\sqrt{a^2-1}}{\sqrt{a^2+2a\cos\theta+1}}+ \frac{\sqrt{a^2+2a\cos\theta+1}}{\sqrt{a^2-1}}          \right )\nonumber\\
&=\frac{a(a+\cos\theta)}{\sqrt{a^2-1}\sqrt{a^2+2a\cos\theta+1}}.
\label{rhs}
\end{align}
Comparison of  (\ref{lhs}) and (\ref{rhs}) shows $\cosh x=e^{-t}\cos y$ and therefore establishes the equivalence of (\ref{pc2}) with the paperclip solution (see Sec. \ref{intro}).

\subsection{The hairclip}{\label{sec:3d}}

In this section, it is shown that the hairclip solution corresponds to the choice $0<a(t)<1$ in the parametric representation (\ref{pc1}) with $g(t)$ chosen to ensure that the solution curve is centred about the the imaginary axis, giving
\begin{equation}
z=-\log\left ( a(t)+e^{i\theta}\right )+\frac{1}{2}\log(1-a^2).
\label{hc1}
\end{equation}
Proceeding as in Sec. \ref{sec:3c} and substituting (\ref{hc1}) into (\ref{schwarz}) yields an ordinary differential equation for $a(t)$ with solution
\begin{equation}
a(t)^2=\frac{1}{1+e^{2t}},
\label{hcasol}
\end{equation}
where an unimportant choice of arbitrary integration constant has been made. Note that $-\infty < t< \infty$ and the solution (\ref{hc1}) and (\ref{hcasol}) can be regarded as an {\it eternal} solution of the curve shortening problem. Again, by considering real and imaginary parts of (\ref{hc1}) and following similar steps as in Sec. \ref{sec:3c}, it can be shown that (\ref{hc1}) and (\ref{hcasol})  is equivalent to the hairclip solution $e^{-t}\cos y=-\sinh x$.

\section{Revisiting the grim reaper solution}
\label{sec:4}
An alternative solution procedure based on solving (\ref{schwarz}) in a steadily moving frame of reference gives the grim reaper solution in the explicit form ${\bar z}=S(z,t)$. Let the $Z$-frame move in the positive real direction with unit speed so that $Z=z-t$. Hence ${\bar Z}={\bar z}-t=f(Z)=f(z-t)$, where $f(Z)$ is the (stationary) Schwarz function of the solution curve in the moving frame. Since ${\bar z}=S(z,t)$, then $S(z,t)=f(z-t)+t$ and the partial derivatives of $S(z,t)$ with respect to $t$ and $z$ are 
\begin{align}
S_{t}&=-f_Z+1\nonumber\\
S_z&=f_Z,\quad {\rm and}\quad S_{zz}=f_{ZZ}.
\label{moving}
\end{align}
In terms of $f(Z)$, (\ref{schwarz}) becomes the ODE
\begin{equation}
f_{ZZ}=f_Z-f_{Z}^2,
\label{ode}
\end{equation}
which has general solution 
\begin{equation}
f(Z)=\log(1+Ke^Z)+C,
\label{odesol1}
\end{equation}
where $K$ and $C$ are constants. Without loss of generality the condition $f(0)=0$ is imposed, which implies the solution curve passes through $Z=0$, and gives $C=-\log(1+K)$. Further, note that as $z\to\infty$, $f\to Z+\log(K/(1+K))$, and so in this limit the solution curve behaves like $y\to(i/2)\log(K/(1+K))$. Thus $|K/(1+K)|=1$. Demanding that the solution is symmetric about the real $Z$-axis implies that $K/(1+K)=-1$ (i.e. $y\to\pm \pi/2$ as $z\to\infty$) and so $K=-1/2$. Finally
\begin{equation}
f(Z)=\log(2-e^Z),
\label{odesol2}
\end{equation}
and the explicit Schwarz function representation of the grim reaper solution is
\begin{equation}
{\bar z}=S(z,t)=\log(2e^t-e^z).
\label{grexpl}
\end{equation}
By taking real and imaginary parts of (\ref{grexpl}) it is straightforward to verify that it is equivalent to the standard grim reaper expression $x=t-\log(\cos y)$.

While it is known that the grim reaper is the only steadily translating solution to the curve shortening problem (e.g. \cite{hall}), the demonstration in this section suggests that the search for steadily translating solutions to more general curve evolution problems in which $v_n=F(\kappa)$, where $F$ is a differentiable function, might be fruitfully pursued by this approach. That is, by formulating the ODE version of $v_n=F(\kappa)$ in a moving frame, and solving the ODE to find explicit solutions of the Schwarz function. 

To give an example of this suppose $F(\kappa)=\kappa^{1/3}$, so that the curve evolution $v_n=\kappa^{1/3}$ gives the Schwarz function PDE 
\begin{equation}
\frac{i}{2}\frac{S_{t}}{(S_z)^{1/2}}=\left [ \frac{i}{2}\frac{S_{zz}}{(S_z)^{3/2}}\right ]^{1/3}.
\label{variant}
\end{equation}
Seeking a solution to (\ref{variant}) translating with unit speed in the positive real direction implies that $S(z,t)=f(Z)+t$ where, from (\ref{variant}),
 \begin{equation}
f_{ZZ}=-\frac{1}{4}(1-f_Z)^3,
\label{variantode}
\end{equation}
where $Z=z-t$. The general solution of (\ref{variantode}) is $f=Z+\sqrt{4A-8Z}+B$, where $A$ and $B$ are constants. Since $f-Z=-2i{\rm Im}(Z)$, squaring and considering the imaginary part gives $B=-2$. Further, requiring the solution curve passes through $Z=0$ gives $A=1$. Thus, in the $z$-plane the solution curve has Schwarz function
 \begin{equation}
S(z,t)=z+\sqrt{4-8(z-t)},
\label{variantsol}
\end{equation}
which upon taking real and imaginary parts implies the curve is the translating parabola $y^2=2(x-t)$. It is straightforward to verify that this translating parabolic solution $\pmb\gamma$ satisfies (\ref{variant}) by directly calculating its curvature and showing $\kappa^{1/3}=1/\sqrt{1+y^2}$, and showing this is the same as the normal velocity $v_n=(\partial {\pmb\gamma}/\partial t).{\textbf{ n}}$.

\section{Remarks}
\label{sec:5}
The Schwarz function formulation of the curve shortening problem in the complex plane, equation (\ref{schwarz}), together with assumptions on the symmetry of the solutions, enables a relatively simple derivation of the fundamental solutions of curve shortening. 
The four solutions can be summarised by the relation
\begin{equation}
z=-\log(a+\zeta)+\frac{1}{2}\log|a^2-1|,
\label{summ}
\end{equation}
where $\zeta=\exp(i\theta)$ and with the following realisations:
\begin{enumerate}
\item $a^2=1/(1-\exp(2t))$, $-\infty <t<0$, $-\pi\le \theta\le \pi$ is the paperclip solution;

\item $a^2=1/(1+\exp(2t))$, $-\infty <t,\theta <\infty$ is the hairclip solution;

\item The limit $t\to -\infty$ of (\ref{summ}) gives $a^2\to 1+\exp(2t)$ and recovers the grim reaper solution (\ref{grparam}). That is, the ancient time limit of the hairclip or paperclip consists of either an ensemble, or pair, of grim reapers which approach each other as $t$ increases forming the hair-- and paperclip solutions respectively.

\item The limit $t\to 0$ of the paperclip solution, gives $a\to\infty$ and $z\to \zeta/a$; that is, a circle with vanishingly small radius. 

\end{enumerate}

Note that solutions of the form (\ref{summ})  with $|a(t)|\le 1$ also arise in the analysis of finger evolution resulting from the instability of an interface separating fluids with different viscosities in the absence of surface tension e.g. \cite{howi2,mineev2}, where superposition of such solutions was used to examine their stability and to construct new solutions for evolving fingers. Owing to the nonlinearity of (\ref{schwarz}) considering a superposition of solutions does not seem to be a way of generating further exact solutions, but it may offer a useful approach to the numerical study of curve shortening of, say, periodic interfaces.

Finally, as noted in  Sec. \ref{sec:2},  curve shortening flow  satisfies, in terms of the Schwarz function, the heat equation 
$S_t=2S_{ss}$ (\ref{heat}). This form of the evolution equation is strongly connected to the geometric evolution equation used by, e.g. \cite{tsai}, to find the fundamental curve shortening solutions. It would be of interest to further pursue the connection between the Schwarz function approach of the present work to   geometric heat equation based methods.

\vspace{1cm}

\noindent I am grateful to Sam Harris for useful discussion related to this work.

%
%



\end{document}